# Mode-Shape Expansion Using Physics-Constrained Gaussian Process Regression

Farid Ghahari[1,2*]

[1] *California Geological Survey, Sacramento, CA 95814, USA*

[2] *The B. John Garrick Institute for the Risk Sciences, University of California, Los Angeles, CA 90095, USA*

## ABSTRACT

This paper addresses the challenge of reconstructing full-field structural mode shapes from sparse sensor data. While Gaussian Process Regression (GPR) offers a robust non-parametric framework for spatial interpolation and uncertainty quantification, standard formulations often yield physically inconsistent mode-shape reconstructions under sparse sensing conditions. A Physics-Constrained Single-Output Gaussian Process (CONS-SOGP) framework is derived that utilizes independent modal kernels while coupling the optimization via a mass-orthogonality penalty. The paper presents derivations for the marginal likelihood, hyperparameter gradients, and penalty coupling. Numerical verification on a multi-degree-of-freedom structure demonstrates that the proposed method overcomes existing limitations in GP-based prediction, providing more accurate and reliable expanded mode shapes.

## INTRODUCTION

Modal identification of structural systems has been a subject of significant interest for many years [1,2]. Modal properties, such as natural frequencies and mode shapes, are essential for tracking temporary or permanent changes in structural systems resulting from gradual aging, extreme events like earthquakes, soil-structure interaction, or environmental factors [3-7]. While several techniques exist to identify modes using a limited number of sensors [8-10], the identified mode shapes remain restricted to the sensor locations. Expanding such sparse representations to a larger (unmeasured) set of degrees of freedom (DOFs) allows for better visualization and denser comparisons—for instance, with analytical mode shapes during Finite Element (FE) model updating. Furthermore, these expanded shapes can be utilized for damage assessment, active or passive control, and response reconstruction (e.g., [11]).

Traditionally, mode shape expansion methods are classified into two groups. The first group relies on the availability of a Finite Element (FE) model or a specific parametric representation—such as cubic polynomials in spline interpolation—to represent the mode shapes [12,13]. However, the reliability of these methods is highly sensitive to modeling uncertainties. The second group consists of data-driven methods, such as projection techniques [14]. Currently, the number of model-based approaches far exceeds that of data-driven methods, as the expansion problem is inherently underdetermined. To account for modeling uncertainty in model-based techniques, a penalty factor representing error in the measured mode shape is commonly employed [15].

Gaussian Process Regression (GPR) [16,17] is a non-parametric, probabilistic model that has been widely adopted across various fields, including earthquake and structural engineering [18-21]. Beyond providing mean predictions, GPR offers a quantified prediction uncertainty, which serves as a critical metric when working with real-world data. While this method has been utilized in several studies for mode shape reconstruction [22,23], current literature—to the author's knowledge—treats modes independently via Single-Output Gaussian Process (SOGP) models. This ignores the physical reality that modes are eigenvectors of the same system matrix. In this paper, these SOGP models are coupled through a mass-orthogonality condition, resulting in the proposed Constrained SOGP (CONS-SOGP) method. Numerical examples demonstrate that this new model provides superior mean estimation with reduced uncertainty. Furthermore, the hyperparameter identifiability problem, which typically manifests as plateaus in the log-marginal likelihood, is resolved by this penalized GPR approach.

[*] E-mail: farid.ghahari@conservation.ca.gov

# SINGLE-OUTPUT GAUSSIAN PROCESS (SOGP)

Assume that the $j$-th mode shape, $f_j(x)$, over the normalized height of the building $x \in [0,1]$ follows a zero-mean[†] Gaussian Process

$$f_j(x) \sim \mathcal{GP}\left(0, k_j(x, x'; \boldsymbol{\vartheta}_j)\right), \tag{1}$$

where $k_j(x, x'; \boldsymbol{\vartheta}_j) = \mathbb{E}[f_j(x) f_j(x')]$ is the prior covariance kernel, with $\mathbb{E}[.]$ denoting the expected value. $\boldsymbol{\vartheta}_j$ is the vector of hyperparameters defining the covariance kernel for the $j$-th mode shape, which will be discussed later. According to the definition of a GP, the vector comprising the mode shape at $n_d$ instrumented floors $\overline{\boldsymbol{x}} = [\bar{x}_1 \quad \bar{x}_2 \quad \cdots \quad \bar{x}_{n_d}]^T$, and the mode shape at any missing floor $x^*$ follows a multi-variate Gaussian distribution given by

$$\begin{bmatrix} f_j(\overline{\boldsymbol{x}}) \\ f_j(x^*) \end{bmatrix} \sim N\left(\begin{bmatrix} \mathbf{0} \\ 0 \end{bmatrix}, \begin{matrix} \mathbf{C}_j(\overline{\boldsymbol{x}}, \overline{\boldsymbol{x}}'; \boldsymbol{\vartheta}_j) + \hat{\sigma}_j^2 \mathbf{I}_{n_d \times n_d} & \boldsymbol{k}_j(\overline{\boldsymbol{x}}, x^*; \boldsymbol{\vartheta}_j) \\ \boldsymbol{k}_j(\overline{\boldsymbol{x}}, x^*; \boldsymbol{\vartheta}_j)^T & k_j(x^*, x^*; \boldsymbol{\vartheta}_j) \end{matrix}\right), \tag{2}$$

where the covariance matrix $\mathbf{C}_j(\overline{\boldsymbol{x}}, \overline{\boldsymbol{x}}'; \boldsymbol{\vartheta}_j)$ is calculated using the kernel function as

$$\mathbf{C}_j(\overline{\boldsymbol{x}}, \overline{\boldsymbol{x}}'; \boldsymbol{\vartheta}_j) = \begin{bmatrix} k_j(\bar{x}_1, \bar{x}_1; \boldsymbol{\vartheta}_j) & k_j(\bar{x}_1, \bar{x}_2; \boldsymbol{\vartheta}_j) & \cdots & k_j(\bar{x}_1, \bar{x}_{n_d}; \boldsymbol{\vartheta}_j) \\ k_j(\bar{x}_2, \bar{x}_1; \boldsymbol{\vartheta}_j) & k_j(\bar{x}_2, \bar{x}_2; \boldsymbol{\vartheta}_j) & \vdots & k_j(\bar{x}_2, \bar{x}_{n_d}; \boldsymbol{\vartheta}_j) \\ \vdots & \cdots & \ddots & \vdots \\ k_j(\bar{x}_{n_d}, \bar{x}_1; \boldsymbol{\vartheta}_j) & k_j(\bar{x}_{n_d}, \bar{x}_2; \boldsymbol{\vartheta}_j) & \cdots & k_j(\bar{x}_{n_d}, \bar{x}_{n_d}; \boldsymbol{\vartheta}_j) \end{bmatrix}, \tag{3}$$

and the vector $\boldsymbol{k}_j(\overline{\boldsymbol{x}}, x^*; \boldsymbol{\vartheta})$ represents the correlation between the measured and predicted points, i.e.,

$$\boldsymbol{k}_j(\overline{\boldsymbol{x}}, x^*; \boldsymbol{\vartheta}_j) = \begin{bmatrix} k_j(\bar{x}_1, x^*; \boldsymbol{\vartheta}_j) \\ k_j(\bar{x}_2, x^*; \boldsymbol{\vartheta}_j) \\ \vdots \\ k_j(\bar{x}_{n_d}, x^*; \boldsymbol{\vartheta}_j) \end{bmatrix}. \tag{4}$$

The term $\hat{\sigma}_j^2$ in Eq. (2) represents the uncertainty (variance) associated with the identified/measured $j$-th mode shape, which is considered to be a zero-mean stationary white Gaussian random process. Using Eq. (2), the conditional distribution of $f_j(x^*)$ given the measurements is a Gaussian distribution as

$$f_j(x^*) | \overline{\boldsymbol{x}}, \overline{\boldsymbol{y}} \sim \mathcal{N}\left(m_{f_j(x^*)}, \sigma^2_{f_j(x^*)}\right), \tag{5}$$

where

$$m_{f_j(x^*)} = \boldsymbol{k}_j^T (\mathbf{C}_j + \hat{\sigma}_j^2 \mathbf{I})^{-1} \overline{\boldsymbol{y}}_j, \tag{6}$$

$$\sigma^2_{f_j(x^*)} = k_j(x^*, x^*; \boldsymbol{\vartheta}_j) - \boldsymbol{k}_j^T (\mathbf{C}_j + \hat{\sigma}_j^2 \mathbf{I})^{-1} \boldsymbol{k}_j, \tag{7}$$

are the posterior mean and variance of the predicted $j$-th mode shape at the normalized height $x^*$. In Eq. (6), $\overline{\boldsymbol{y}}_j = [\varphi_{1,j} \quad \varphi_{2,j} \quad \cdots \quad \varphi_{n_d,j}]^T$ is the vector of $j$-th mode shape measured at $n_d$ instrumented floors. Provided that the kernel function and its hyperparameters for each mode are known, Eqs. (6) and (7) can be used to expand each mode shape using identified/measured elements at the instrumented floors.

The purpose of this study is to improve the GP prediction by connecting the mode shapes through a physical constraint instead of treating them independently. Therefore, the optimal choice of the kernel function is beyond the scope of this study. In this study, a commonly used kernel function called the Squared Exponential (SE) kernel is used for all mode shapes. This function is $C^\infty$ smooth and defined as

$$k_j\left(x, x'; \boldsymbol{\vartheta}_j = [\gamma_j \quad \beta_j]^T\right) = \gamma_j^2 \, e^{-\frac{(x-x')^2}{2\beta_j^2}}, \tag{8}$$

where $\gamma_j$ and $\beta_j$ are hyperparameters of the kernel function for the $j$-th mode shapes, representing the signal standard deviation and correlation length, respectively. To find these hyperparameters within this Single-Output Gaussian Process (SOGP) framework, the Negative Log-Marginal Likelihood (NLML) of the data is minimized. Given $n_d$ instrumented floor and the assumed measurement/identification noise, the marginal likelihood is

$$p(\overline{\mathbf{y}}_j|\overline{\mathbf{x}}, \boldsymbol{\vartheta}_j) = \mathcal{N}(\overline{\mathbf{y}}_j|\mathbf{0}, \mathbf{C}_j + \hat{\sigma}_j^2\mathbf{I}). \tag{9}$$

---

[†] Assuming a zero mean is common in GP modeling, as trends in the data can be removed a priori.

Therefore, the NLML to be minimized is

$$\mathcal{L}_j(\boldsymbol{\vartheta}_j) = \frac{1}{2}\overline{\mathbf{y}}_j^T(\mathbf{C}_j + \hat{\sigma}_j^2\mathbf{I})^{-1}\overline{\mathbf{y}}_j + \frac{1}{2}\log|\mathbf{C}_j + \hat{\sigma}_j^2\mathbf{I}| + \frac{n_d}{2}\log 2\pi. \tag{10}$$

Optimization via quasi-Newton methods requires the partial derivatives of $\mathcal{L}_j(\boldsymbol{\vartheta}_j)$ with respect to each hyperparameter $\theta_{i,j}$ (where $\theta_{1,j} = \gamma_j$ *and* $\theta_{2,j} = \beta_j$), which is

$$\frac{\partial \mathcal{L}_j(\boldsymbol{\vartheta}_j)}{\partial \theta_{i,j}} = -\frac{1}{2}\mathrm{tr}\left(\left[\boldsymbol{\alpha}_j\boldsymbol{\alpha}_j^T - (\mathbf{C}_j + \hat{\sigma}_j^2\mathbf{I})^{-1}\right]\frac{\partial \mathbf{C}_j}{\partial \theta_{i,j}}\right). \tag{11}$$

where $\boldsymbol{\alpha}_j = (\mathbf{C}_j + \hat{\sigma}_j^2\mathbf{I})^{-1}\overline{\mathbf{y}}_j$ and

$$\frac{\partial \mathbf{C}_j}{\partial \gamma_j} = \frac{2}{\gamma_j}\mathbf{C}_j, \tag{12}$$

$$\frac{\partial \mathbf{C}_j}{\partial \beta_j} = \mathbf{C}_j \odot \frac{(\overline{\mathbf{x}}\mathbf{1}_{n_d}^T - \mathbf{1}_{n_d}\overline{\mathbf{x}}^T)^2}{\beta_j^3}. \tag{13}$$

The symbol $\odot$ in Eq. (13) represents element-wise multiplication.

# PHYSICS-CONSTRAINED SOGP (CONS-SOGP)

Complete mode shapes are theoretically linearly independent; however, their independent reconstruction from sparse data often does not respect this constraint. Herein, the orthogonality of expanded mode shapes with respect to the mass matrix $\mathbf{M}$, which is usually easy to estimate for a building structure, is added as a constraint to the expansion process, and the estimation is carried out in an interconnected approach.

Assume there are $n_m$ modes to be expanded. Similar to the previous section, $n_m$ independent kernels are considered, but their hyperparameters —collected in the vector $\boldsymbol{\Theta} = [\boldsymbol{\vartheta}_1^T, \ldots, \boldsymbol{\vartheta}_{n_m}^T]^T$—are estimated through a global objective function defined as the sum of all $n_m$ NLMLs and a penalty function:

$$\mathcal{J}(\boldsymbol{\Theta}) = \lambda\,\mathcal{P}(\boldsymbol{\Phi}) + \sum_{j=1}^{n_m}\mathcal{L}_j(\boldsymbol{\vartheta}_j), \tag{14}$$

where $\boldsymbol{\Phi} = [\boldsymbol{\mu}_1 \quad \cdots \quad \boldsymbol{\mu}_{n_m}]$ is the collection of the posterior mean estimations of mode shapes obtained from Eq. (6) using the hyperparameters derived through solving this new optimization problem. The penalty $\mathcal{P}$ is the squared Frobenius norm of the mass orthogonality violation:

$$\mathcal{P}(\boldsymbol{\Phi}) = \left\|\boldsymbol{\Phi}^T\mathbf{M}\boldsymbol{\Phi} - \mathbf{I}_{n_m\times n_m}\right\|_F^2. \tag{15}$$

The gradient of the penalty with respect to a specific mode shape $\boldsymbol{\mu}_j$ is:

$$\frac{\partial \mathcal{P}}{\partial \boldsymbol{\mu}_j} = 4\mathbf{M}\left(\sum_{i=1}^{n_m}\boldsymbol{\mu}_i(\boldsymbol{\mu}_i^T\mathbf{M}\boldsymbol{\mu}_j) - \boldsymbol{\mu}_j\right), \tag{16}$$

which can be collectively written as follows for all mode shapes

$$\frac{\partial \mathcal{P}}{\partial \boldsymbol{\Phi}} = 4\mathbf{M}\boldsymbol{\Phi}(\boldsymbol{\Phi}^T\mathbf{M}\boldsymbol{\Phi} - \mathbf{I}). \tag{17}$$

Therefore, the derivative of the penalty function with respect to the hyperparameters can be calculated using the chain rule as

$$\frac{\partial \mathcal{P}}{\partial \theta_{i,j}} = \mathrm{tr}\left[\left(\frac{\partial \mathcal{P}}{\partial \boldsymbol{\Phi}}\right)^T\frac{\partial \boldsymbol{\Phi}}{\partial \theta_{i,j}}\right], \tag{18}$$

in which

$$\frac{\partial \boldsymbol{\Phi}}{\partial \theta_{i,j}} = \left[\mathbf{0} \quad \cdots \quad \mathbf{0} \quad \frac{\partial \boldsymbol{\mu}_j}{\partial \theta_{i,j}} \quad \mathbf{0} \quad \cdots \quad \mathbf{0}\right]. \tag{19}$$

According to Eq. (6), $\boldsymbol{\mu}_j = \mathbf{K}_j{}^T\left(\mathbf{C}_j + \hat{\sigma}_j^2\mathbf{I}\right)^{-1}\overline{\boldsymbol{y}}_j$ where $\mathbf{K}_j = \mathbf{K}_j\left(\overline{\boldsymbol{x}}, \boldsymbol{x}^*; \boldsymbol{\vartheta}_j\right)$ is now an $n_d \times N$ matrix, and $N$ is the total number of floors listed in vector $\boldsymbol{x}^*$. Because the SE kernel function is used in this study, it is straightforward to show that,

$$\frac{\partial \boldsymbol{\mu}_j}{\partial \gamma_j} = \frac{\mathbf{2}}{\gamma_j}\mathbf{K}_j{}^T\left(\mathbf{C}_j + \hat{\sigma}_j^2\mathbf{I}\right)^{-1}\hat{\sigma}_j^2\mathbf{I}\,\boldsymbol{\alpha}_j, \tag{20}$$

$$\frac{\partial \boldsymbol{\mu}_j}{\partial \beta_j} = \frac{\mathbf{1}}{\beta_j}\left[\frac{\partial \mathbf{K}_j{}^T}{\partial \ln \beta_j}\boldsymbol{\alpha}_j - \mathbf{K}_j{}^T\left(\mathbf{C}_j + \hat{\sigma}_j^2\mathbf{I}\right)^{-1}\frac{\partial \mathbf{C}_j}{\partial \ln \beta_j}\boldsymbol{\alpha}_j\right], \tag{21}$$

where

$$\frac{\partial \mathbf{K}_j{}^T}{\partial \ln \beta_j} = \mathbf{K}_j{}^T \odot \frac{\left(\boldsymbol{x}^*\mathbf{1}_{n_d}^T - \mathbf{1}_N\overline{\boldsymbol{x}}^T\right)^2}{\beta_j^2}, \tag{22}$$

$$\frac{\partial \mathbf{C}_j}{\partial \ln \beta_j} = \mathbf{C}_j \odot \frac{\left(\overline{\boldsymbol{x}}\mathbf{1}_{n_d}^T - \mathbf{1}_{n_d}\overline{\boldsymbol{x}}^T\right)^2}{\beta_j^2}. \tag{23}$$

After finding the optimal hyperparameters of all modes by minimizing Eq. (14), the mean expanded mode shapes and their uncertainties are calculated for each mode using Eq. (6) and (7).

The penalty function introduced in Eq. (15) is only valid when the mode shapes are mass-normalized. However, mode shapes identified from vibration data are not inherently mass-normalized. Furthermore, because these mode shapes are sparse, it is not possible to make them mass-normalized prior to expansion. To resolve this issue, the expanded mode shape means, $\boldsymbol{\mu}_j$s, are mass normalized as $\widetilde{\boldsymbol{\mu}}_j = \boldsymbol{\mu}_j / \sqrt{\boldsymbol{\mu}_j^T\mathbf{M}\boldsymbol{\mu}_j}$ during the optimization, and then the penalty function is calculated. As long as the mass-normalized mode shapes $\widetilde{\boldsymbol{\Phi}} = [\widetilde{\boldsymbol{\mu}}_1 \quad \cdots \quad \widetilde{\boldsymbol{\mu}}_{n_m}]$ are used in Eq. (17), the only adjustment required for the the gradient of the penalty function is the modification of the term $\partial\boldsymbol{\Phi}/\,\partial\theta_{i,j}$. Utilizing the chain rule, $\frac{\partial \widetilde{\boldsymbol{\mu}}_j}{\partial \theta_{i,j}} = \frac{\partial \widetilde{\boldsymbol{\mu}}_j}{\partial \boldsymbol{\mu}_j}\frac{\partial \boldsymbol{\mu}_j}{\partial \theta_{i,j}}$ where $\frac{\partial \widetilde{\boldsymbol{\mu}}_j}{\partial \boldsymbol{\mu}_j}$ is calculated as

$$\frac{\partial \widetilde{\boldsymbol{\mu}}_j}{\partial \boldsymbol{\mu}_j} = \frac{\mathbf{1}}{\sqrt{\boldsymbol{\mu}_j^T\mathbf{M}\boldsymbol{\mu}_j}}\left(\mathbf{I} - \frac{\boldsymbol{\mu}_j\boldsymbol{\mu}_j^T\mathbf{M}}{\boldsymbol{\mu}_j^T\mathbf{M}\boldsymbol{\mu}_j}\right). \tag{24}$$

Finally, it is recommended to use the logarithms of the hyperparameters when minimizing the negative log-likelihood to ensure these parameters remain positive. Additionally, since mode shapes are usually zero at the base (in the absence of soil-structure interaction effects), a virtual observation of zero can be added to the list of instrumented floors. To prevent numerical instabilities, it is suggested to use a small jitter value instead of absolute zero for the base noise. MATLAB [24] implementations of both SOGP and CONS-SOGP for mode shape expansion are provided in the Appendix for the reader's convenience.

# NUMERICAL EXAMPLE

To verify the proposed CONS-SOGP method for mode shape expansion and demonstrate its superiority over the standard SOGP, a shear building model with 53 floors above the base is numerically modeled. This model is constructed to have natural periods close to the 52-story building that has been extensively studied by the author [25-27]. To make the case more challenging, a sudden 80% reduction in stiffness is considered at the 15$^{\text{th}}$ floor, representing a damaged building. Ground truth mode shapes are calculated using an eigen analysis of $\mathbf{M}^{-1}\mathbf{K}$, where $\mathbf{M}$ and $\mathbf{K}$ are the mass and stiffness matrices, respectively. To ensure the mode shapes are not initially mass-normalized, each complete mode shape is normalized to its norm. It is assumed the first 5 modes are only identified/available at 10$^{\text{th}}$, 20$^{\text{th}}$, 30$^{\text{th}}$, 40$^{\text{th}}$, and 53$^{\text{rd}}$ floors (yielding 6 physical measurement levels). Random noise with a standard deviation equal to 2% of the norm of each mode shape vector is added to the modal deformations at these specific floors.

To carry out the expansion analysis using both the SOGP and CONS-SOGP methods, $\gamma_j$ and $\beta_j$ are initially set to 1 and 0.1, respectively. The parameter $\gamma_j$ is constrained between 0.1 and 10, while $\beta_j$ is limited between 0.02 and 1.5. These are relatively wide ranges considering that both the mode shapes and floor levels are normalized for the analysis. **Figure 1** shows the comparison between the expanded mode shapes (mean and 95% confidence interval) utilizing the SOGP method and the ground truth mode shapes. As seen in this figure, the expanded mode shapes for the first two modes represent the true mode shapes very well, and the associated uncertainty is very small. In the remaining three modes, however, the GP model exhibits spiky behavior, where the predicted mean collapses toward zero between measurement points. Consequently, the uncertainty is very large between sensors. This phenomenon is common in GP applications and is caused by the flatness of the NLML. To support this observation, the NLML is plotted across the entire range of the correlation length $\beta_j$, assuming the signal variance $\gamma_j^2$ is set at a constant value for each mode (these values correspond to the optimal values obtained from the CONS-SOGP method). As observed in **Figure 2**, a global minimum exists in the NLML function for the first mode and slightly for the second mode; however, for the other three modes, there is a wide flat region for small values of $\beta_j$. Therefore, the minimization algorithm sets $\beta_j$ to its lower limit, which causes the observed spiky behavior in the predicted means. As mentioned earlier, this is a well-documented issue, and several solutions have been proposed to address it. The most common approach is to penalize the NLML function with an additional term (e.g., the inverse of the correlation length) (see, e.g., [28]). However, adding such a purely mathematical penalty term to the NLML is subjective and lacks physical theoretical support.

The CONS-SOGP method proposed in this study penalizes the NLML function using a physics-based constraint: mode shape orthogonality with respect to the mass matrix. As expected, this penalty not only mitigates the observed spiky behavior of the GP model, but it also significantly improves the mean prediction by incorporating physically meaningful constraints into the model. **Figure 3** illustrates the mode shape expansions generated through CONS-SOGP. As shown, the mean predictions align closely with the actual mode shapes. Furthermore, the associated uncertainties are reasonable and adequately cover the discrepancy between the estimated and true mode shapes wherever the predicted mean deviates from the true mode shape (e.g., at the damaged 15$^{th}$ floor). The optimal hyperparameters identified through each of these two methods are reported in Table 1. As anticipated, the correlation length remains well above the lower limits when utilizing the proposed CONS-SOGP method.

In the example presented here, the penalty factor λ is set to 1000, but it should be adjusted for each case specifically, depending on the number of mode shapes and measurement points. A rule of thumb is to adjust the penalty factor such that both terms in Eq. (14) have approximately the same value. This can be achieved through a few iterations.

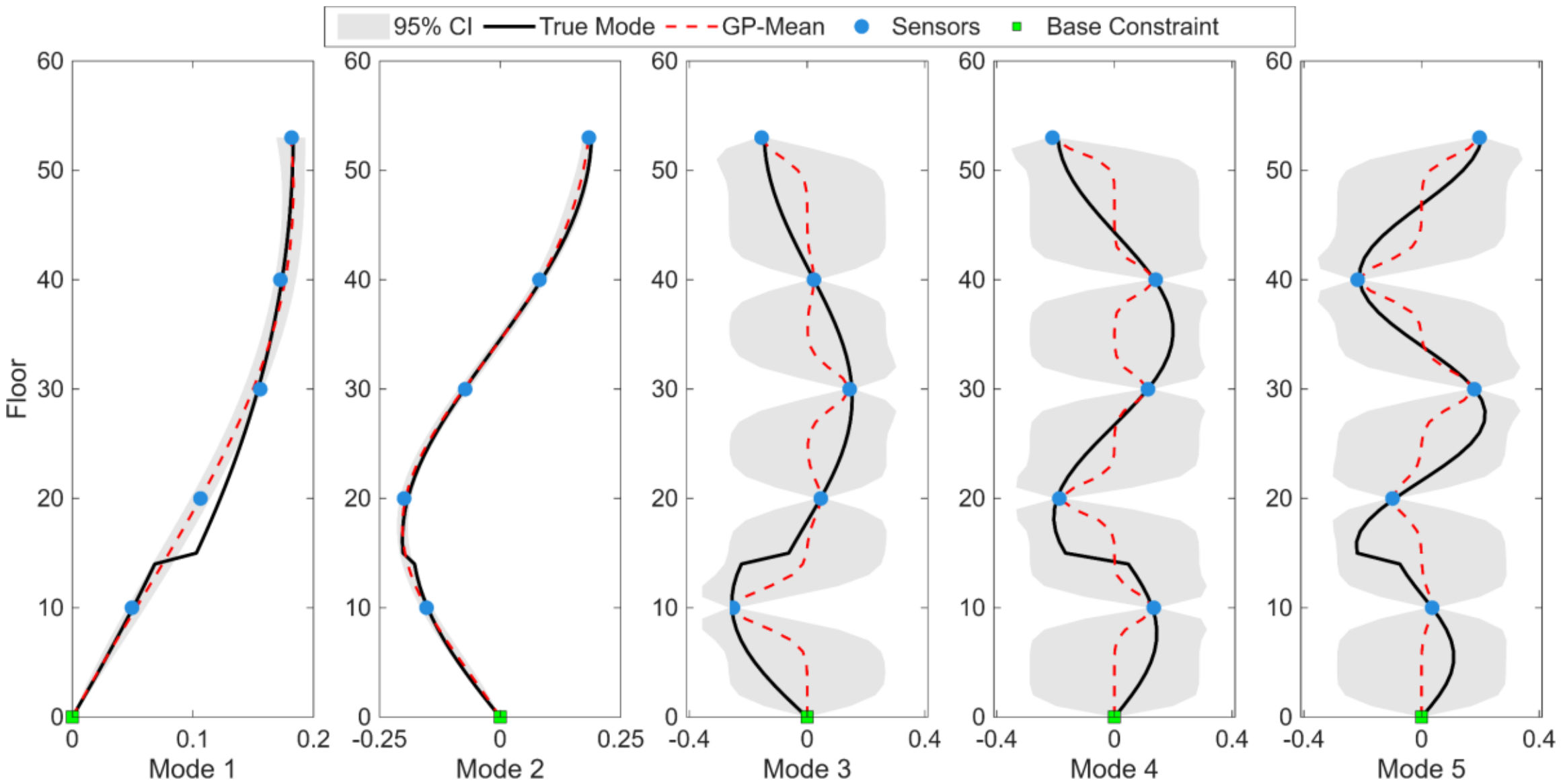

**Figure 1.** Comparison between ground-truth and SOGP-expanded mode shapes.

**Table 1.** Identified hyperparameters for the squared exponential (SE) kernel function.

| Mode | SOGP | | CONS-SOGP | |
|---|---|---|---|---|
| | $\gamma_j$ | $\beta_j$ | $\gamma_j$ | $\beta_j$ |
| 1 | 2.27 | 0.77 | 2.27 | 0.42 |
| 2 | 0.99 | 0.33 | 0.92 | 0.27 |
| 3 | 0.84 | 0.03 | 0.88 | 0.11 |
| 4 | 0.82 | 0.02 | 2.78 | 0.23 |
| 5 | 0.82 | 0.03 | 1.03 | 0.11 |

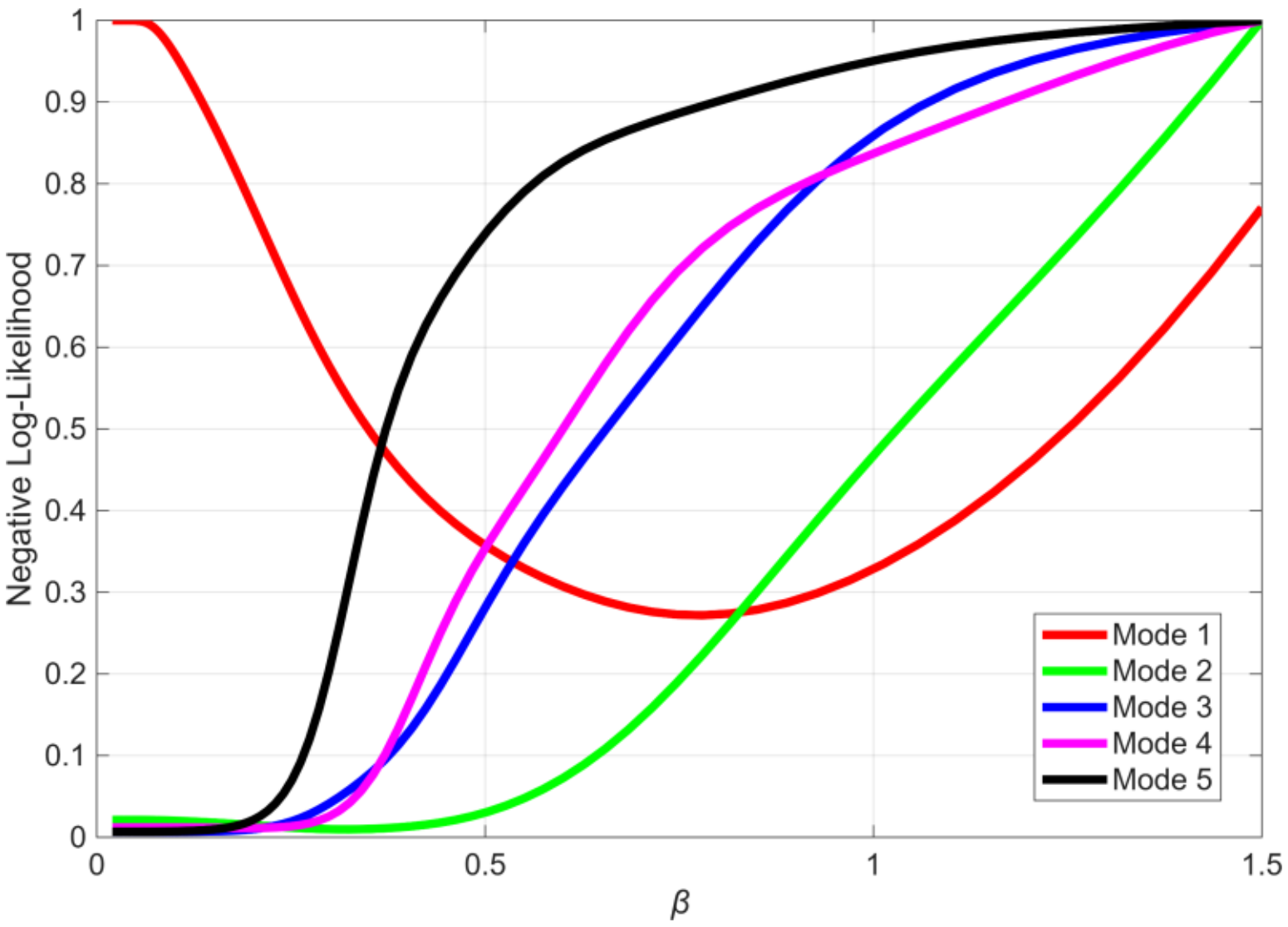

**Figure 2.** Variation of the negative log-marginal likelihood (NLML) in the SOGP method versus the kernel function correlation length for a fixed signal variance.

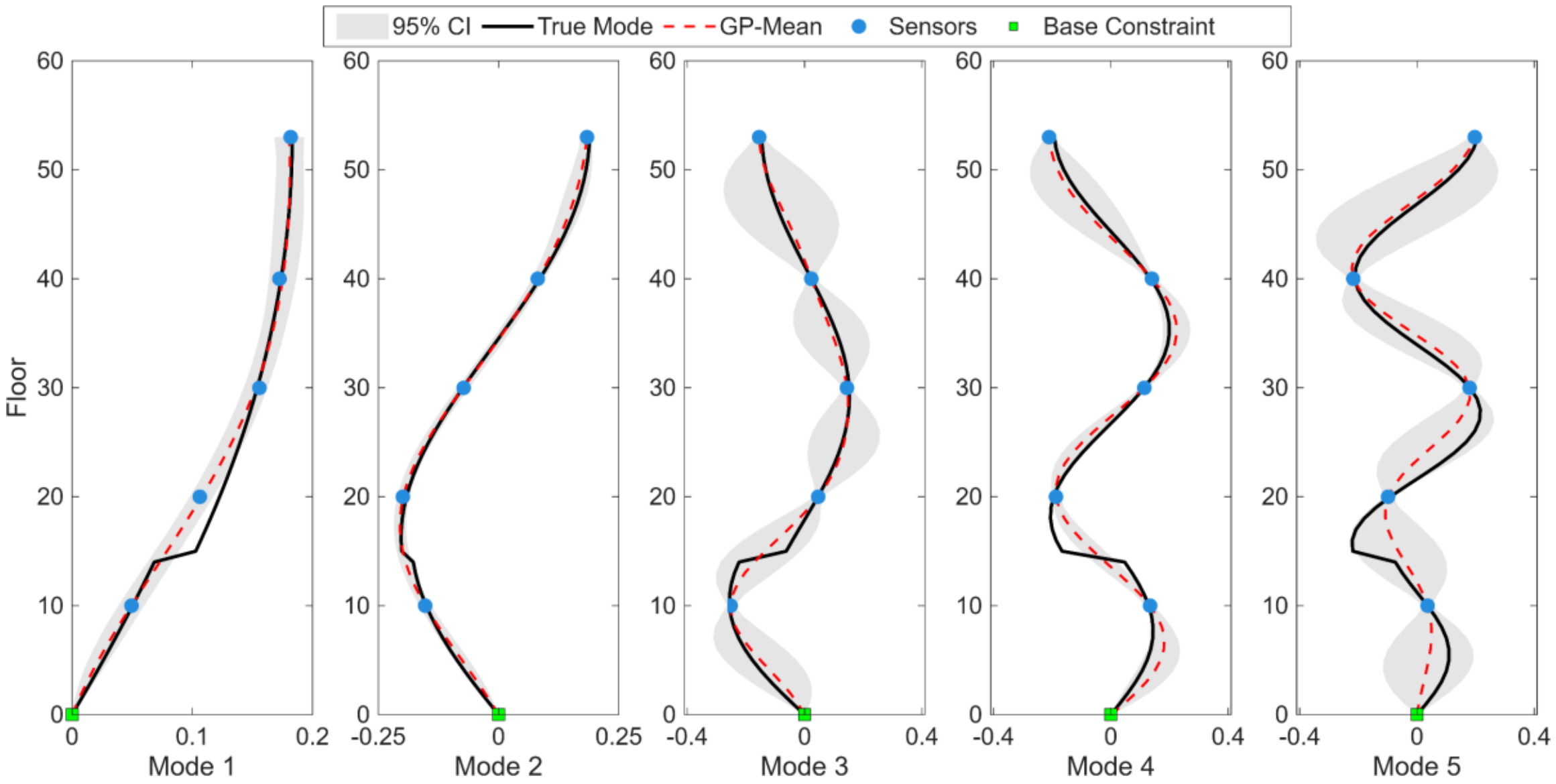

**Figure 3.** Comparison between ground-truth and CONS-SOGP-expanded mode shapes.

# CONCLUSIONS

This paper proposes a probabilistic method for mode shape expansion. While the Single-Output Gaussian Process (SOGP) has already been used for mode shape expansion, it treats each mode separately even though they are eigenvectors of the same system matrix. Furthermore, the standard SOGP suffers from likelihood flatness, which can result in spiky mean estimations. To address these two limitations, the method proposed here connects all modes through the fact that mode shapes are orthogonal with respect to the mass matrix, which is usually available or can be accurately estimated. The orthogonality violation is added as a penalty function to the negative log-likelihood to improve the estimation of the kernel function hyperparameters for all modes. A simulated example demonstrates the superiority of the proposed method, termed Constrained SOGP (CONS-SOGP).

# ACKNOWLEDGMENTS

This study was supported by the California Geological Survey (CGS), which is gratefully acknowledged. Any opinions, findings, and conclusions expressed in this material are solely those of the author and do not reflect the views of the California Department of Conservation.

# APPENDIX

```
function [Phi_Mean,Phi_Std,Theta] = SOGP_func(Y_phys,obs_floors,N,noise_std_phys)
	% Input:
	% Y_phys : Measured mode shapes at instrumented degrees of freedom (n_d x n_m)
	% obs_floors : Vector of instrumented floor indices (n_d x 1)
	% N : Total number of floors (excluding the base)
	% noise_std_phys : Standard deviation of the noise in physical units  (scalar or 1 x n_m)

	% Output:
	% Phi_Mean : Mean estimate of the expanded mode shapes ((N+1) x n_m)
	% Phi_Std : Standard deviation of the expanded mode shapes ((N+1) x n_m)
	% Theta : Optimum hyperparameters (n_m x 2)

	% Constants:
	base_noise = 1e-6;
	theta0 = [log(1.0), log(0.1)]; % gamma, beta
	lb = [log(0.1), log(0.02)];
	ub = [log(10.0), log(1.5)];

	[n_d,n_m] = size(Y_phys);
	Y_aug = [zeros(1, n_m); Y_phys];
	x_norm = (0:N)' / N;
	x_obs_aug = x_norm([0; obs_floors(:)] + 1);
	Phi_Mean = zeros(N+1, n_m);
	Phi_Std = zeros(N+1, n_m);
	if isscalar(noise_std_phys)
		noise_std_phys = repmat(noise_std_phys, 1, n_m);
	end
	opts = optimoptions('fmincon', 'Display', 'off', 'Algorithm', 'sqp',...
	'SpecifyObjectiveGradient', true,'OptimalityTolerance', 1e-7,'StepTolerance', 1e-7);
	Theta = [];
	for k = 1:n_m
		y_k = Y_aug(:, k);
		scale_factor = std(Y_phys(:, k));
		y_norm_in = y_k / scale_factor;
		sn_phys_norm = noise_std_phys / scale_factor;
		sn_vec = [base_noise; ones(n_d, 1) * sn_phys_norm];
		[theta_opt, ~] = fmincon(@(t) nlml(t, x_obs_aug, y_norm_in, sn_vec), theta0, [], [], [], [], lb, ub, [], opts);
		gamma = exp(theta_opt(1));
		beta = exp(theta_opt(2));
		dist2_full = (x_norm - x_norm').^2;
		dist2_obs = (x_obs_aug - x_obs_aug').^2;
		dist2_fo = (x_norm - x_obs_aug').^2;
		K_full = gamma^2 * exp(-0.5 * dist2_full / beta^2);
		K_obs = gamma^2 * exp(-0.5 * dist2_obs / beta^2);
		K_fo = gamma^2 * exp(-0.5 * dist2_fo / beta^2);
		Ky = K_obs + diag(sn_vec.^2);
		L = chol(Ky, 'lower');
		alpha = L' \ (L \ y_norm_in);
		mu_norm = K_fo * alpha;
		v = L \ K_fo';
		var_norm = diag(K_full - v' * v);
		std_norm = sqrt(max(0, var_norm));
		Phi_Mean(:, k) = mu_norm * scale_factor;
		Phi_Std(:, k) = std_norm * scale_factor;
		Theta = [Theta;[gamma,beta]];
	 end
end

%%%%%%%%%%%%%%%%%%%%%%%%%%%%%%%%%%%%%%%%%%%%%%%%%%%%%%%%%%%%%%%%%%%%
function [L_val,Grad] = nlml(theta, x_obs, y, sn_vec)
	gamma = exp(theta(1));
	beta = exp(theta(2));
	n = length(x_obs);
	x_diff = x_obs - x_obs';
	sq_dist = x_diff.^2;
	K = gamma^2 * exp(-0.5 * sq_dist / beta^2);
	Ky = K + diag(sn_vec.^2);
	[L_mat, p] = chol(Ky, 'lower');
	if p > 0
		L_val = 1e12;
		Grad = zeros(size(theta));
		return;
	end
	alpha = L_mat' \ (L_mat \ y);
	L_val = 0.5 * (y' * alpha) + sum(log(diag(L_mat))) + (n/2) * log(2*pi);
	if nargout > 1
		invL = L_mat \ eye(n);
		invKy = invL' * invL;
		W = alpha * alpha' - invKy;
		dK_dlogsf = 2 * K;
		dK_dlogl = K .* (sq_dist / beta^2);
		Grad = zeros(2, 1);
		Grad(1) = -0.5 * sum(sum(W .* dK_dlogsf'));
		Grad(2) = -0.5 * sum(sum(W .* dK_dlogl'));
	end
end
```

**Figure A1.** MATLAB implementation of the SOGP method for mode shape expansion.

```matlab
function [Phi_Mean, Phi_Std, Theta] = CSOGP_func(Y_phys, obs_floors, noise_std_phys, M, lambda)
    % Input:
    % Y_phys : Measured mode shapes at instrumented degrees of freedom (n_d x n_m)
    % obs_floors : Vector of instrumented floor indices (n_d x 1)
    % noise_std_phys : Standard deviation of the noise in physical units (scalar or 1 x n_m)
    % M: : Mass matrix (N x N)
    % lambda : Penalty factor

    % Output:
    % Phi_Mean : Mean estimate of the expanded mode shapes ((N+1) x n_m)
    % Phi_Std : Standard deviation of the expanded mode shapes ((N+1) x n_m)
    % Theta : Optimum hyperparameters (n_m x 2)

    % Constants
    theta0 = [log(1.0), log(0.1)]; % [gamma, beta]
    lb = [log(0.1), log(0.02)];
    ub = [log(10.0), log(1.5)];
    base_noise = 1e-6;
    [N,~] = size(M);
    [n_d, n_m] = size(Y_phys);
    x_norm = (0:N)' / N;
    x_obs_aug = x_norm([0; obs_floors(:)] + 1);
    scale_vec = std(Y_phys, 0, 1);
    Y_norm_aug = [zeros(1, n_m); Y_phys] ./ scale_vec;
    if isscalar(noise_std_phys)
        sn_phys_norm = repmat(noise_std_phys, 1, n_m) ./ scale_vec;
    else
        sn_phys_norm = noise_std_phys ./ scale_vec;
    end
    theta0 = repmat(theta0, n_m, 1);
    lb = repmat(lb, n_m, 1); ub = repmat(ub, n_m, 1);
    opts = optimoptions('fmincon', 'Display', 'iter', 'Algorithm', 'sqp', SpecifyObjectiveGradient', true, 'OptimalityTolerance', 1e-10, ...
    'StepTolerance', 1e-10, 'CheckGradients', true, 'FiniteDifferenceType', 'central');
    [theta_opt, ~] = fmincon(@(t) nlml(t, x_obs_aug, x_norm, Y_norm_aug, n_m, M, scale_vec, lambda, sn_phys_norm, base_noise), theta0, [], [], [], [], lb, ub, [], opts);
    Phi_Mean = zeros(N+1, n_m); Phi_Std = zeros(N+1, n_m);
    for k = 1:n_m
        gamma = exp(theta_opt(k, 1));
        beta = exp(theta_opt(k, 2));
        sn_vec = [base_noise; ones(n_d, 1) * sn_phys_norm(k)];
        sq_dist_obs = (x_obs_aug - x_obs_aug').^2;
        sq_dist_full = (x_norm - x_norm').^2;
        sq_dist_fo = (x_norm - x_obs_aug').^2;
        Ky = gamma^2 * exp(-0.5 * sq_dist_obs / beta^2) + diag(sn_vec.^2);
        K_fo = gamma^2 * exp(-0.5 * sq_dist_fo / beta^2);
        K_full = gamma^2 * exp(-0.5 * sq_dist_full / beta^2);
        L = chol(Ky, 'lower');
        alpha = L' \ (L \ Y_norm_aug(:,k));
        Phi_Mean(:,k) = (K_fo * alpha) * scale_vec(k);
        v = L \ K_fo';
        Phi_Std(:,k) = sqrt(max(0, diag(K_full - v'*v))) * scale_vec(k);
    end
    Theta = exp(theta_opt);
end
%%%%%%%%%%%%%%%%%%%%%%%%%%%%%%%%%%%%%%%%%%%%%%%%%%%%%%%%%%%%%%%%%%%%%%%%%%%%%%%%%%%%%%%%%%%%%%%%%%%%%%%%%%%
function [J, Grad] = nlml(theta, x_obs, x_full, Y_aug, n_m, M, scale_vec, lambda, sn_phys, base_noise)
    n_obs = length(x_obs);
    n_full = length(x_full);
    Phi_mu = zeros(n_full, n_m);
    Grad = zeros(size(theta));
    total_NLML = 0;
    sq_dist_obs = (x_obs - x_obs').^2;
    sq_dist_fo = (x_full - x_obs').^2;
    for k = 1:n_m
        gamma = exp(theta(k, 1));
        beta = exp(theta(k, 2));
        sn_vec = [base_noise; ones(n_obs-1, 1) * sn_phys(k)];
        K_obs = gamma^2 * exp(-0.5 * sq_dist_obs / beta^2);
        Ky = K_obs + diag(sn_vec.^2);
        [L, p] = chol(Ky, 'lower');
        if p > 0, J = 1e15; Grad = zeros(size(theta)); return; end
        alpha = L' \ (L \ Y_aug(:,k));
        total_NLML = total_NLML + 0.5*(Y_aug(:,k)'*alpha) + sum(log(diag(L))) + (n_obs/2)*log(2*pi);
        K_fo = gamma^2 * exp(-0.5 * sq_dist_fo / beta^2);
        Phi_mu(:,k) = (K_fo * alpha) * scale_vec(k);
    end
    Phi_act = Phi_mu(2:end, :);
    modal_mass = diag(Phi_act' * M * Phi_act);
    Phi_norm = Phi_act ./ sqrt(modal_mass');
    Orth_Err = Phi_norm' * M * Phi_norm - eye(n_m);
    Penalty = lambda * sum(Orth_Err(:).^2);
    J = total_NLML + Penalty;
    if nargout > 1
        dJ_dPhi_norm = 4 * lambda * (M * Phi_norm * Orth_Err);
        for k = 1:n_m
            gamma = exp(theta(k, 1));
            beta = exp(theta(k, 2));
            sn_vec = [base_noise; ones(n_obs-1, 1) * sn_phys(k)];
            K_obs = gamma^2 * exp(-0.5 * sq_dist_obs / beta^2);
            Ky = K_obs + diag(sn_vec.^2);
            L = chol(Ky, 'lower');
            alpha = L' \ (L \ Y_aug(:,k));
            K_fo_act = gamma^2 * exp(-0.5 * sq_dist_fo(2:end, :) / beta^2);
            invKy = solve_inv(L);
            W = alpha * alpha' - invKy;
            Grad(k, 1) = -0.5 * sum(sum(W .* (2*K_obs)'));
            Grad(k, 2) = -0.5 * sum(sum(W .* (K_obs .* (sq_dist_obs / beta^2))'));
            dPhiact_dloggamma = 2 * (K_fo_act * (L' \ (L \ (diag(sn_vec.^2) * alpha)))) * scale_vec(k);
            dKfo_dlogbeta = K_fo_act .* (sq_dist_fo(2:end,:) / beta^2);
            dKy_dlogbeta = K_obs .* (sq_dist_obs / beta^2);
            term1 = dKfo_dlogbeta * alpha; term2 = K_fo_act * (L' \ (L \ (dKy_dlogbeta * alpha)));
            dPhiact_dlogbeta = (term1 - term2) * scale_vec(k);
            dPhinorm_dloggamma = 1 / sqrt(modal_mass(k)) * (eye(n_full-1)-1/modal_mass(k)*Phi_act(:,k) * Phi_act(:,k)' * M) * dPhiact_dloggamma;
            dPhinorm_dlogbeta = 1 / sqrt(modal_mass(k)) * (eye(n_full-1)-1/modal_mass(k)*Phi_act(:,k) * Phi_act(:,k)' * M) * dPhiact_dlogbeta;
            Grad(k, 1) = Grad(k, 1) + dJ_dPhi_norm(:,k)' * dPhinorm_dloggamma;
            Grad(k, 2) = Grad(k, 2) + dJ_dPhi_norm(:,k)' * dPhinorm_dlogbeta;
        end
    end
end
%%%%%%%%%%%%%%%%%%%%%%%%%%%%%%%%%%%%%%%%%%%%%%%%%%%%%%%%%%%%%%%%%%%%%%%%%%%%%%%%%%%%%%%%%%%%%%%%%%%%%%%%%%%
function invA = solve_inv(L)
        n = size(L, 1); invL = L \ eye(n); invA = invL' * invL;
end
```

**Figure A2.** MATLAB implementation of the CONS-SOGP method for mode shape expansion.